\documentclass[12pt,a4paper]{article}
\usepackage{pifont}
\usepackage{amsfonts}
\usepackage{mathrsfs}
\usepackage{amssymb}
\usepackage{amsmath}

\newtheorem{theorem}{Theorem}[section]

\newtheorem{lemma}[theorem]{Lemma}

\makeatletter

\newcommand{\Rmnum}[1]{\expandafter\@slowromancap\romannumeral #1@}
\makeatother

\begin{document}
\title{Embedding dendriform dialgebra into its universal enveloping
Rota-Baxter algebra\footnote{Supported by the NNSF of China (Nos.
10771077; 10911120389) and the NSF of Guangdong Province
 (No. 06025062).}}
\author{
 Yuqun
Chen \ and   Qiuhui Mo \\
{\small \ School of Mathematical Sciences, South China Normal
University}\\
{\small Guangzhou 510631, P.R. China}\\
{\small yqchen@scnu.edu.cn}\\
{\small scnuhuashimomo@126.com}}

\date{}

\maketitle \noindent\textbf{Abstract:}    In this paper, by using
Gr\"{o}bner-Shirshov bases for Rota-Baxter algebras, we prove that
every dendriform dialgebra over a field of characteristic 0 can be
embedded into its universal enveloping Rota-Baxter algebra of weight
0.



\section{Introduction}

 Let $F$ be a field. A dendriform dialgebra (see
\cite{Lo}) is a $F$-module $D$ with two binary operations $\prec$
and $\succ$ such that for any $x,y,z\in D$,
\begin{eqnarray}\label{e1}
\nonumber&&(x\prec y)\prec z=x\prec(y\prec z+y\succ z)\\
&&(x\succ y)\prec z=x\succ(y\prec z)\\
\nonumber&&(x\prec y+x\succ y)\succ z=x\succ(y\succ z)
\end{eqnarray}

Let $A$ be an associative algebra over $F$ and $\lambda\in F$. Let a
$F$-linear operator $P: A\rightarrow A$ satisfy
\begin{eqnarray}\label{e2}
P(x)P(y)=P(P(x)y)+P(xP(y))+\lambda P(xy), \forall x,y\in A.
\end{eqnarray}
Then $A$ is called a Rota-Baxter algebra of weight $\lambda$.

The free Rota-Baxter algebra generated by a nonempty set $X$,
denoted by $RB(X)$, and the free dendriform algebras  generated by
$X$, denoted by $D(X)$, are given by K. Ebrahimi-Fard and L. Guo
\cite{EG08a}.

Suppose that $(D, \prec,\succ)$ is a dendriform dialgebra over $F$
with a linear basis $X=\{x_i|i\in I\}$. Let $x_i\prec x_j=\{x_i\prec
x_j\}, x_i\succ x_j=\{x_i\succ x_j\}$, where $\{x_i\prec x_j\}$ and
$\{x_i\succ x_j\}$ are linear combinations of $x\in X$. Then $D$ has
an expression by generator and defining relations
$$
D=D(X|x_i\prec x_j=\{x_i\prec x_j\}, x_i\succ x_j=\{x_i\succ x_j\},
x_i,x_j\in X).
$$

Denote by
$$
U(D)=RB(X|x_iP(x_j)+\lambda x_ix_j=\{x_i\prec x_j\},
P(x_i)x_j=\{x_i\succ x_j\}, x_i,x_j\in X).
$$
Then $U(D)$ is the universal enveloping Rota-Baxter algebra of $D$,
where $\lambda \in F$, see \cite{EG08a}.

\ \

The study of Rota-Baxter algebra originated from the probability
study of Glenn Baxter in 1960 and was developed further by Cartier
and the school of Rota in the 1960s and 1970s. This structure
appeared also in the Lie algebra context as the operator form of the
classical Yang-Baxter equation started in the 1980s. Since then,
Rota-Baxter algebra has experienced a quite remarkable renascence
and found important theoretical developments and applications in
mathematical physics, operads, number theory and combinatorics, see,
for example, \cite{AGKO, Baxter, Cartier, E02, EG08b, Rota}.

The dendriform dialgebra  was introduced by J.-L. Loday \cite{Lo} in
1995 with motivation from algebraic K-theory, and was further
studied in connection with several areas in mathematics and physics,
including operads, homology, Hopf algebras, Lie and Leibniz
algebras, combinatorics, arithmetic and quantum field theory, see
\cite{EG08a,Lo04}.

In the theory of Lie algebras, the Poincare-Birkhoff-Witt theorem
(Poincare \cite{Poincare}, Birkhoff \cite{Birkhoff}, Witt
\cite{Witt}, frequently contracted to PBW theorem) is a fundamental
result giving an explicit description of the universal enveloping
algebra of a Lie algebra. The term ``PBW type theorem" or even ``PBW
theorem" may also refer to various analogues of the original
theorem, see, for example, \cite{Cohn,Higgins}.

I.P. Shestakov \cite{Shest} proved that an Akivis algebra can be
embedded into its universal enveloping non-associative algebra. M.
Aymon and P.-P. Grivel \cite{A} proved that a Leibniz algebra can be
embedded into its universal enveloping dialgebra. P.S. Kolesnikov
\cite{K} proved that every (finite dimensional) Leibniz algebra can
be embedded into current conformal algebra over the algebra of
linear transformations of a (finite dimensional) linear space. As a
corollary,  a new proof of the theorem on injective embedding of a
Leibniz algebra into an associative dialgebra is obtained and, more
explicitly,  an analogue of the PBW theorem for Leibniz algebras in
 \cite{Lo}.

Motivated by the above, in this paper, we prove the following
theorem.

\begin{theorem}\label{t1}
Every dendriform dialgebra over a field of characteristic 0 can be
embedded into its universal enveloping Rota-Baxter algebra of weight
0. In another words, such a dendriform dialgebra is isomorphic to a
dendriform subdialgebra  of a Rota-Baxter algebra of weight 0.

\end{theorem}

Composition-Diamond lemma for Rota-Baxter algebras is established by
L.A. Bokut, Yuqun Chen and Xueming Deng in a recent paper
\cite{BCD}. In this paper, by using this lemma, we prove the above
theorem.

\section{Preliminaries}

In this section, we introduce some notations which are related
Gr\"{o}bner-Shirshov bases for Rota-Baxter algebras, see \cite{BCD}.

Let $X$ be a nonempty set, $S(X)$ the free semigroup generated by
$X$ without identity and $P$ a symbol of a unary operation. For any
two nonempty sets $Y$ and $Z$, denote by
$$
{\Lambda}_{P}{(Y,Z)}=(\cup_{r\geq 0}(YP(Z))^{r}Y)\cup(\cup_{r\geq
1}(YP(Z))^{r})\cup(\cup_{r\geq 0}(P(Z)Y)^{r}P(Z))\cup(\cup_{r\geq
1}(P(Z)Y)^{r}),
$$
where for a set $T$,  $T^0$ means the empty set.

Define
\begin{eqnarray*}
\Phi_{0}&=&S(X)\\
\vdots\ \ & &\ \ \ \ \vdots\\
\Phi_{n}&=&{\Lambda}_{P}(\Phi_{0},\Phi_{n-1}) \\
\vdots\ \ & &\ \ \ \ \vdots\\
\end{eqnarray*}
Then
$$
\Phi_{0} \subset \cdots \subset\Phi_{n} \subset
 \cdots
$$
Let
$$
\Phi(X)=\cup_{n\geq 0} \Phi_{n}.
$$
Clearly, $P(\Phi(X))\subset \Phi(X)$. If $u\in X\cup P({\Phi}(X))$,
then $u$ is called prime.  For any $u\in \Phi(X)$, $u$ has a unique
form
 $u=u_{1}u_{2}\cdots u_{n}$ where $u_{i}$  is prime, $i=1, 2,\dots,n$, and
 $u_i,u_{i+1}$ can not both have forms as $P(u'_i)$ and
 $P(u'_{i+1})$.

For any $u\in \Phi(X)$ and for a set $T\subseteq X\cup \{P\}$,
denote by $deg_{T}(u)$ the number of occurrences of $t\in T$ in $u$.
 Let
$$
Deg(u)=(deg_{\{P\}\cup X}(u), deg_{\{P\}}(u)).
$$
We order $Deg(u)$ lexicographically.

In the following, we always assume that $F$ is a field of
characteristic 0.

Let $F\Phi(X)$ be a free  $F$-module with $F$-basis $\Phi(X)$ and
$\lambda \in F$ a fixed element. Extend linearly $P: \
F\Phi(X)\rightarrow F\Phi(X),\ u\mapsto P(u)$ where $u\in \Phi(X)$.

Now we define the multiplication in $F\Phi(X)$.

Firstly, for $u,v\in X\cup P(\Phi(X))$, define

$u\cdot v= \left\{ \begin{array}
              {l@{\quad}l}
              P(P(u')\cdot v')+P(u'\cdot P(v'))+\lambda P(u'\cdot v'), & \mbox{if}\ u=P(u'),v=P(v'); \\
              uv, & \mbox{otherwise}.
              \end{array} \right. $

Secondly, for any $u=u_{1}u_{2}\cdots u_{s}, v=v_{1}v_{2}\cdots
v_{l}\in \Phi(X)$ where $u_{i}, v_{j}$ are prime, $i=1,2,\dots,s,
j=1,2,\dots,l$, define
$$
u\cdot v=u_{1}u_{2}\cdots u_{s-1}(u_{s}\cdot v_{1})v_{2}\cdots
v_{l}.
$$

Equipping with the above concepts, $F\Phi(X)$ is the free
Rota-Baxter algebra with weight $\lambda$ generated by $X$, see
\cite{EG08a}.

We denote by $RB(X)$ the free Rota-Baxter algebra with weight
$\lambda$ generated by $X$.

Let  $\cal N^+$ be the set of positive integers.

Let the notations be as before. We have to order $\Phi(X)$. Let $X$
be a well-ordered set. Let us define an ordering $>$ on $\Phi(X)$ by
induction on the $Deg$-function.

For any $u,v\in \Phi(X)$, if $Deg(u)>Deg(v)$, then $u>v$.

If $Deg(u)=Deg(v)=(n,m)$, then we define $u>v$ by induction on
$(n,m)$.

If $(n,m)=(1,0),$ then $u,v\in X$ and we use the ordering on $X$.
Suppose that for $(n,m)$ the ordering is defined where
$(n,m)\geq(1,0)$. Let $(n,m)<(n',m')=Deg(u)= Deg(v)$.  If $u,v \in
P(\Phi(X))$, say $u=P(u')$ and $v=P(v')$, then $u>v$ if and only if
$u'>v'$ by induction. Otherwise $u=u_1u_2\cdots u_l$ and
$v=v_1v_2\cdots v_s$ where $l>1$ or $s>1$, then $u>v$ if and only if
$(u_1, u_2, \dots, u_l)>(v_1, v_2, \dots, v_s )$ lexicographically
by induction.

It is clear that $>$ is a well ordering on $\Phi(X)$, see
\cite{BCD}. Throughout this paper, we will use this ordering.

Let $\star$ be a symbol and $\star\notin X$. By a
$\star$-Rota-Baxter word we mean any expression in $\Phi(X\cup
\{\star\})$ with only one occurrence of $\star$. The set of all
$\star$-Rota-Baxter words on $X$ is denoted by $\Phi^\star (X)$.

Let $u$ be a $\star$-Rota-Baxter word and $s\in RB(X)$. Then we call
$$
u|_{s}=u|_{\star\mapsto s}
$$
an $s$-Rota-Baxter word. For short, we  call $u|_{s}$ an $s$-word.

Note that the ordering $>$ is monomial in the sense that for any
$u,v\in \Phi(X),\ w\in\Phi^{\star}(X)$,
$$
u>v\Longrightarrow \overline{w|_{u}}>\overline{w|_{v}}
$$ where $w|_u=w|_{\star\mapsto u}$ and $w|_v=w|_{\star\mapsto v}$, see \cite{BCD}, Lemma 3.4.

If $\overline{u|_{s}}=u|_{\overline s}$, then we call $u|_{s}$ a
normal $s$-word.

Now, for any $0\neq f\in RB(X)$, $f$ has the leading term $\bar{f}$
and $f=\alpha_{1}\bar{f}+\sum_{i=2}^{n}\alpha_iu_i$ where
$\bar{f},u_i\in \Phi(X), \bar{f}>u_i, 0\neq\alpha_1, \alpha_i \in
F$. Denote by $lc(f)$ the coefficient of the leading term $\bar{f}$.
If $lc(f)=1$, we call $f$ monic.

Let $f,g \in RB(X)$ be monic with $\overline{f}=u_{1}u_{2}\cdots
u_{n}$ where each $u_{i} $ is prime. Then, there are four kinds of
compositions.
\begin{description}
\item{(i)} If $u_{n}\in P(\Phi(X))$, then we define composition of
right multiplication  as\ \ $f\cdot u$ \ \ where $u\in P(\Phi(X)).$

\item{(ii)} If $u_{1}\in P(\Phi(X))$, then we define composition of
left multiplication as \ \ $u\cdot f$ \ \  where \ $u\in
P(\Phi(X)).$

\item{(iii)} If there exits a $w=\overline{f}a=b
\overline{g}$ \ where $fa$ is normal $f$-word and $bg$ is normal
$g$-word, $a, b \in \Phi(X)$ and $deg_{\{P\}\cup
X}(w)<deg_{\{P\}\cup X}(\overline{f})+deg_{\{P\}\cup
X}(\overline{g})$, then we define the intersection composition of
$f$ and $g$ with respect to $w$ as $(f,g)_{w}=f\cdot a-b\cdot g.$

\item{(iv)} If there exists a $w=\overline{f}=u|_{\overline{g}}$ where
$u\in \Phi^{\star}(X)$, then we define the inclusion composition of
$f$ and $g$ with respect to $w$ as $(f,g)_{w}=f-u|_{g}.$
\end{description}

We call $w$ in $(f,g)_{w}$ the ambiguity with respect to $f$ and
$g$.

Let $S\subset RB(X)$ be a set of monic polynomials. Then the
composition $h$ is called trivial modulo $(S,w)$, denoted by
$h\equiv 0 \ \ mod (S,w)$, if
$$
h=\sum_{i}\alpha_iu_i|_{s_{i}}
$$
where each $\alpha_i\in F$,  $s_i\in S$, $u_i|_{s_i}$ is normal
$s_i$-word and $u_i|_{\overline{s_i}}<\bar h$, and $\bar h=w$ if $h$
is a composition of left (right) multiplication.

In general, for any two polynomials $p$ and $q$, $ p\equiv q \ \ mod
(S,w) $ means that $ p-q\equiv0 \ \ mod (S,w)$.

$S$ is called a Gr\"{o}bner-Shirshov basis  in  $RB(X)$ if each
composition  is trivial  modulo $S$ and responding $w$.

\begin{theorem}\label{CDL}(\cite{BCD}, Composition-Diamond lemma for Rota-Baxter
algebras) Let $RB(X)$ be a free Rota-Baxter algebra over a field of
characteristic 0 and $S$  a set of monic polynomials in $RB(X)$, $
>$ the monomial ordering on $\Phi(X)$ defined as before and $Id(S)$ the Rota-Baxter ideal of
$RB(X)$ generated by $S$. Then the following statements are
equivalent.
\begin{enumerate}
\item[(I)] $S $ is a Gr\"{o}bner-Shirshov basis in $RB(X)$.
\item[(II)] $  f\in Id(S)\Rightarrow
\bar{f}=u|_{\overline{s}}$  for some $u \in \Phi^{\star}(X)$,\ $s\in
S$.
\item[(III)] $Irr(S) = \{ u\in \Phi(X) |  u \neq
v|_{\bar{s}}, s\in S, v|_{s}\ \mbox{is normal}\  s\mbox{-word}\}$ is
a $F$-basis of $RB(X|S)\\=RB(X)/Id(S)$.
\end{enumerate}
\end{theorem}

If a subset $S$ of $RB(X)$ is not a Gr\"{o}bner-Shirshov basis then
one can add all nontrivial compositions of polynomials of $S$ to
$S$. Continuing this process repeatedly, we finally obtain a
Gr\"{o}bner-Shirshov basis $S^{comp}$ that contains $S$. Such a
process is called Shirshov algorithm.

\section{The proof of Theorem \ref{t1}}
In this section, we assume that $RB(X)$ is the free Rota-Baxter
algebra generated by $X=\{x_i\mid i\in I\}$ with weight $\lambda=0$.

\begin{lemma}\label{l0}
For any $u,v\in \Phi(X)$, we have
$\overline{P(u)P(v)}=max\{\overline{P(P(u)v)},\overline{P(uP(v))}\}$.
\end{lemma}
{\bf Proof:} By Rota-Baxter formula (\ref{e2}), we may assume that
$P(P(u)v)=\sum n_iu_i,\ P(uP(v))=\sum m_jv_j$, where $n_i,m_j\in
{\cal N^+}, \ u_i,v_j\in \Phi(X)$. Since the characteristic of $F$
is 0, the result follows. $ \hfill \blacksquare $

 \ \

Denote by
\begin{eqnarray*}
&&F_{1}=\{x_iP(x_j)-\{x_i\prec x_j\}\ | \ i,j\in I\},\\
&&F_{2}=\{P(x_i)x_j-\{x_i\succ x_j\}\ | \ i,j\in I\}.\\
&&Irr(F_{1}\cup F_{2})) = \{ u\in \Phi(X) | u \neq v|_{\bar{s}},
s\in F_{1}\cup F_{2}, v|_{s}\ \mbox{is normal}\ s\mbox{-word}\},\\
&&\Phi_{1}(X)=\Phi(X)\cap Irr(F_{1}\cup F_{2}).
\end{eqnarray*}
For a polynomial $f=\Sigma_{i=1}^n\alpha_iu_i\in RB(X)$, where each
$0\neq\alpha_i\in F,\ u_i\in \Phi(X)$, denote the set $\{u_i, \
1\leq i\leq n\}$ by $supp(f)$.

\begin{lemma}\label{l2}
\begin{description}
\item 1)\ Let $f=P(x_i)u$, $g=vP(x_j)$, where $i,j\in I, \  u,v\in
\Phi_{1}(X)\setminus X$. Then $f\equiv \Sigma\alpha_iu_i
 \ mod(F_{1}\cup F_{2}, \overline{f})$ and $g\equiv \Sigma\beta_iv_i
 \ mod(F_{1}\cup F_{2}, \overline{g})$, where for any $i,\ \alpha_i,\beta_i\in F,\ u_i, v_i\in \Phi_{1}(X)\setminus X$.

\item 2) \ Let $f=P(u)P(v)$, $g=P(v')P(u')$, where $u, u'\in
\Phi_{1}(X)\setminus X,\ v, v'\in \Phi_{1}(X)$. Then $f\equiv
\Sigma\alpha_iP(u_i)
 \ mod(F_{1}\cup F_{2}, \overline{f})$ and $g\equiv \Sigma\beta_iP(v_i)
 \ mod(F_{1}\cup F_{2}, \overline{g})$, where for any $i,\ \alpha_i,\beta_i\in F,\ u_i, v_i\in \Phi_{1}(X)\setminus
 X$.
\end{description}
\end{lemma}

\noindent\textbf{Proof.} \textbf{1)} \ We prove only the case
$f=P(x_i)u\equiv \Sigma\alpha_iu_i \ mod(F_{1}\cup F_{2},
\overline{f})$. Another case is similar.

We use induction on $n=deg_{\{P\}\cup X}(u)$. Since $u\in
\Phi_{1}(X)\setminus X$, we have $n\geq 2$.

Assume that $n=2$. Then either $u=x_jx_k$ or $u=P(x)$, $x_j, x_k,
x\in X$. If $u=x_jx_k$, we have $f=P(x_i)x_jx_k\equiv\{x_i\succ
x_j\}x_k \ mod(F_{1}\cup F_{2}, \overline{f})$, and $supp(\{x_i\succ
x_j\}x_k)\subset \Phi_{1}(X)\setminus X$. If $u=P(x)$, we have
$f=P(x_i)P(x)=P(P(x_i)x)+P(x_iP(x))\equiv P(\{x_i\succ
x\})+P(\{x_i\prec x\})\ mod(F_{1}\cup F_{2}, \overline{f})$, and
$supp(P(\{x_i\succ x\})), supp(P(\{x_i\prec x\}))\subset
\Phi_{1}(X)\setminus X$.

For $n>2$,   there are three cases to consider.
\begin{description}

\item (I) $u=x_ju_1, \ x_j\in X$. Then there are two subcases to consider.
\begin{description}
\item (i) $u_1=x_ku_2, \ x_k\in X$. Then $f=P(x_i)u\equiv \{x_i\succ
x_j\}x_ku_2\ \ mod(F_{1}\cup F_{2}, \overline{f})$, where
$supp(\{x_i\succ x_j\}x_ku_2)\subset \Phi_{1}(X)\setminus X$.
\item (ii) $u_1=P(v)u_2$. Since $u\in
\Phi_{1}(X)\setminus X$, we get that $v\not\in X$. Thus,
$f=P(x_i)u\equiv \{x_i\succ x_j\}P(v)u_2\ mod(F_{1}\cup F_{2},
\overline{f})$, where $supp(\{x_i\succ x_j\}P(v)u_2)\subset
\Phi_{1}(X)\setminus X$.
\end{description}

\item (II) $u=P(u_1)$. Then $f=P(x_i)u=P(x_i)P(u_1)=P(P(x_i)u_1)+P(x_iP(u_1)
)$. Let $P(x_i)u_1\equiv \sum\gamma_iw_i\ \ mod(F_{1}\cup F_{2},
\overline{P(x_i)u_1})$ and $x_iP(u_1)\equiv \sum\gamma_i'w_i'\ \
mod(F_{1}\cup F_{2}, \overline{x_iP(u_1)})$, where all $w_i,w_i'\in
\Phi_{1}(X)$. By using Lemma \ref{l0},
$\overline{f}=\overline{P(x_i)P(u_1)}\geq
\overline{P(P(x_i)u_1)}=P(\overline{P(x_i)u_1})\geq P(w_j)$
 and similarly, $\overline{f}\geq P(w_j')$ for any $j, j'$. Then
$f\equiv \Sigma \gamma_iP(w_i)+ \Sigma \gamma_i'P(w_i')\
mod(F_{1}\cup F_{2}, \overline{f})$, where $P(w_i),P(w_i')\in
\Phi_{1}(X)\setminus X$.

\item (III) $u=P(u_1)u_2$, where $u_2$ is not empty. Then
$u_1\not\in X$, and $u_2=x_ju_3$ for some $x_j\in X$ since $u\in
\Phi_{1}(X)\setminus X$. Therefore,
$f=P(x_i)u=P(x_i)P(u_1)u_2=P(P(x_i)u_1)u_2+P(x_iP(u_1) )u_2$. For
$P(x_iP(u_1) )u_2$, we have $P(x_iP(u_1) )u_2\in
\Phi_{1}(X)\setminus X$. For $P(P(x_i)u_1)u_2$, since $u_1\not\in
X$, by induction on $n$, we get that $P(x_i)u_1\equiv
\Sigma\gamma_iv_i \ mod(F_{1}\cup F_{2}, \overline{P(x_i)u_1})$,
where $v_i\in \Phi_{1}(X)\setminus X$. By using Lemma \ref{l0},
$\overline{f}=\overline{P(x_i)P(u_1)u_2}\geq
\overline{P(P(x_i)u_1)u_2}=P(\overline{P(x_i)u_1})u_2\geq
P(v_i)u_2$. As a result $P(P(x_i)u_1)u_2\equiv
\Sigma\gamma_iP(v_i)u_2 \ mod(F_{1}\cup F_{2}, \overline{f})$ and
$P(v_i)u_2\in \Phi_{1}(X)\setminus X$.
\end{description}

\textbf{2)}\  We only prove the case $f=P(u)P(v)$. Another case is
similarly to prove.

We use induction on $n=deg_{\{P\}\cup X}(P(u)P(v))$. Since $u\in
\Phi_{1}(X)\setminus X$, we have $n\geq 5$.

Assume that $n=5$. Then either $u=x_ix_j$ and $v=x$ or $u=P(x_i)$
and $ v=x$, where $x_i, x_j, x\in X$.

If $u=x_ix_j$ and $v=x$, we have
$f=P(u)P(v)=P(x_ix_j)P(x)=P(P(x_ix_j)x)+P(x_ix_jP(x))\equiv
P(P(x_ix_j)x)+P(x_i\{x_j\prec x\})\ mod(F_{1}\cup F_{2},
\overline{f})$, and $(\{P(x_ix_j)x\}\cup supp(x_i\{x_j\prec
x\}))\subset \Phi_{1}(X)\setminus X$.

If $u=P(x_i)$ and $v=x$, we have
$f=P(u)P(v)=P(P(x_i))P(x)=P(P(P(x_i))x)+P(P(x_i)P(x))=
P(P(P(x_i))x)+P(P(P(x_i)x))+P(P(x_iP(x)))\equiv
P(P(P(x_i))x)+P(P(\{x_i\succ x\}))+P(P(\{x_i\prec x\}))\
mod(F_{1}\cup F_{2}, \overline{f})$, and $(\{P(P(x_i))x\}\cup
supp(P(\{x_i\succ x\}))\cup supp(P(\{x_i\prec x\})))\subset
\Phi_{1}(X)\setminus X$.

For $n>5$, since $f=P(u)P(v)=P(uP(v))+P(P(u)v)$ and by Lemma
\ref{l0}, it is sufficient to prove that $P(uP(v))\equiv
\Sigma\alpha_iP(u_i)
 \ mod(F_{1}\cup F_{2}, \overline{P(uP(v))})$, $P(P(u)v)\equiv \Sigma\alpha_iP(v_i)
 \ mod(F_{1}\cup F_{2}, \overline{P(P(u)v)})$, where $u_i, v_i\in \Phi_{1}(X)\setminus
 X$.

For $P(uP(v))$, there are two cases to consider.
\begin{description}

\item (I) $u=u_1x_i, \ x_i\in X$. Then there are two subcases to consider.
\begin{description}

\item (i) $v\not\in X$. Then $P(uP(v))=P(u_1x_iP(v))$ and $u_1x_iP(v)\in \Phi_{1}(X)\setminus
 X$.

\item (ii) $v=x_j\in X$. Then $P(uP(v))=P(u_1x_iP(x_j))\equiv P(u_1\{x_i\prec
x_j\})$. If $u_1=u_2x$ for some $x\in X$, then $P(uP(v))\equiv
P(u_2x\{x_i\prec x_j\})$ where $supp(u_2x\{x_i\prec x_j\})\subset
\Phi_{1}(X)\setminus
 X$. If $u_1=u_2P(u_3)$, then $u=u_1x_i=u_2P(u_3)x_i$ and $u_3\not\in X$. Then $P(uP(v))\equiv
P(u_1\{x_i\prec x_j\})\equiv P(u_2P(u_3)\{x_i\prec x_j\})$ where
$supp(u_2P(u_3)\{x_i\prec x_j\})\subset \Phi_{1}(X)\setminus
 X$.

\end{description}

\item (II) $u=u_1P(u_2)$. Then there are two subcases to consider.
\begin{description}

\item (i) $u_2=x_i\in X$. Since $u\in \Phi_{1}(X)\setminus
 X$, we have  $u=P(x_i)$. As
 a result, $P(uP(v))=P(P(x_i)P(v))$. Since $P(v)\not\in X$, the result follows
 from $1)$.

\item (ii) $u_2\not\in X$. Then $u_2\in \Phi_{1}(X)\setminus
 X$ and $P(uP(v))=P(u_1P(u_2)P(v))$. By
induction on  $n$, $P(u_2)P(v)\equiv \Sigma\alpha_iP(v_i) \ \
mod(F_{1}\cup F_{2}, \overline{P(u_2)P(v)})$, where $v_i\in
\Phi_{1}(X)\setminus
 X$. Then $P(uP(v))=P(u_1P(u_2)P(v))\equiv
 \Sigma\alpha_iP(u_1P(v_i))$ and $u_1P(v_i)\in \Phi_{1}(X)\setminus
 X$.
\end{description}
\end{description}

For $P(P(u)v)$, there are also two cases to consider.
\begin{description}
\item (I) $v=x_iv_1, \ x_i\in X$. Then $P(P(u)v)=P(P(u)x_iv_1)$
and  $P(u)x_iv_1\in \Phi_{1}(X)\setminus
 X$.

\item (II) $v=P(v_1)v_2$. Then $P(P(u)v)=P(P(u)P(v_1)v_2)$ with $v_1\in \Phi_{1}(X)$. By induction on
$n$, we get that $P(u)P(v_1)\equiv \Sigma\alpha_iP(u_i)
 \ mod(F_{1}\cup F_{2}, \overline{P(u)P(v_1)})$, where $u_i\in \Phi_{1}(X)\setminus
 X$. Then $P(P(u)v)\equiv
 \Sigma\alpha_iP(P(u_i)v_2)$ and $P(u_i)v_2 \in \Phi_{1}(X)\setminus
 X$.
\end{description}
The proof is completed.
 $ \hfill \blacksquare $

\ \

\begin{lemma}\label{t2}
Let  $S=F_1\cup F_2\cup F_3$, where
\begin{eqnarray*}
F_{3}&=&\{u_{0}P(v_1)u_{1}\cdots P(v_{n})u_{n}\ |\ u_{0},u_{n}\in
X^*,\ u_{i}\in X^*\backslash
\{1\}, 1\leq i<n,\ \\
&&\ \ v_{j}\in \Phi_1(X)\setminus X,\
 1\leq j\leq n; \  |u_{0}|\geq 2\ \mbox{ if }\ n=0 \},
\end{eqnarray*}
where for any $u\in X^*, \ |u|$ is the length of $u$, $X^*$ is the
free monoid generated by $X$.

Then $S$ is a Gr\"{o}bner-Shirshov basis in $RB(X)$.
\end{lemma}

\noindent\textbf{Proof.} \ The ambiguities of all possible
compositions of the polynomials  in $S$ are only as below:
\begin{enumerate}
\item[$f_1\wedge f_2$]\ \ $f_1\in F_1, \ f_2\in F_2$, and $w=x_iP(x_j)x_k,\ i,j,k\in
I$.

\item[$f_2\wedge f_1$]\ \ $f_1\in F_1, \ f_2\in F_2$, and $w=P(x_i)x_jP(x_k),\ i,j,k\in
I$.

\item[$f_2\wedge f_3$]\ \ $f_2\in F_2, \ f_3\in F_3$, and $w=P(x_i)x_ju_{0}P(v_1)u_{1}\cdots
P(v_{n})u_{n},\ u_{0},u_{n}\in X^*,\ u_{k}\in X^*\backslash\{1\},\
v_{l}\in \Phi_{1}(X)\setminus X,\ i,j\in I,\ n\geq0,\ 1\leq k<n,\
1\leq l\leq n$. When $n=0,\ |u_{0}^{(1)}|\geq 1$.

\item[$f_3\wedge f_1$]\ \  $f_1\in F_1, \ f_3\in F_3$, and $w=u_{0}P(v_1)u_{1}\cdots
P(v_{n})u_{n}x_iP(x_j), \ u_{0},u_{n}\in X^*,\ u_{k}\in
X^*\backslash\{1\},\ v_{l}\in \Phi_{1}(X)\setminus X,\ i,j\in I,\
n\geq0,\ 1\leq k<n,\ 1\leq l\leq n$. When $n=0,\ |u_{0}^{(1)}|\geq
1$.

\item[$f_3\wedge f_3'$]\ \  $f_3, f_3'\in F_3$. There are three ambiguities, one is for the intersection composition and
two are for the inclusion composition.
\end{enumerate}

All possible compositions of left and right multiplication are: $
f_1P(u)$,\  $ P(u)f_2$,\  $ f_3P(u)$ and $P(u)f_3$, where $f_i\in
F_i,\ u\in \Phi(X),\ i=1,2,3$.

Now we prove that all the compositions are trivial.

For $f_1\wedge f_2$,  let $f=x_iP(x_j)-\{x_i\prec x_j\}$,
$g=P(x_j)x_k-\{x_j\succ x_k\},\ i,j,k\in I$. Then $w=x_iP(x_j)x_k$
and
\begin{eqnarray*}
(f,g)_w&=&x_iP(x_j)x_k-\{x_i\prec
x_j\}x_k-(x_iP(x_j)x_k-x_i\{x_j\succ x_k\})\\
&=&x_i\{x_j\succ x_k\}-\{x_i\prec x_j\}x_k \\
&\equiv&0\ \  mod(F_3, w).
\end{eqnarray*}

For $f_2\wedge f_1$, let $f=P(x_i)x_j-\{x_i\succ x_j\}$,
$g=x_jP(x_k)-\{x_j\prec x_k\},\ i,j,k\in I$. Then
$w=P(x_i)x_jP(x_k)$ and by equation (\ref{e1}),
\begin{eqnarray*}
(f,g)_w&=&P(x_i)x_jP(x_k)-\{x_i\succ x_j\}P(x_k)-P(x_i)(x_jP(x_k)-\{x_j\prec x_k\})\\
&=&P(x_i)\{x_j\prec x_k\}-\{x_i\succ x_j\}P(x_k)\\
&\equiv&\{x_i\succ\{x_j\prec x_k\}\}-\{\{x_i\succ x_j\}\prec x_k\}\\
&\equiv&0\ \  mod(S, w).
\end{eqnarray*}

For $f_3\wedge f_1$,  let $f=u_{0}P(v_1)u_{1}\cdots
P(v_{n})u_{n}x_i$, $g=x_iP(x_j)-\{x_i\prec x_j\}$, \ $u_{0},u_{n}\in
X^*, u_{k}\in X^*\backslash\{1\}, v_{l}\in \Phi_{1}(X)\setminus X,
i,j\in I, n\geq0, 1\leq k<n,  1\leq l\leq n$, and $|u_{0}|\geq 1$ if
$n=0$. Then $w=u_{0}P(v_1)u_{1}\cdots P(v_{n})u_{n}x_iP(x_j)$ and
\begin{eqnarray*}
(f,g)_w&=&u_{0}P(v_1)u_{1}\cdots
P(v_{n})u_{n}\{x_i\prec x_j\}\\
&\equiv&0\ \  mod(S, w).
\end{eqnarray*}

For $f_2\wedge f_3$,  the proof is similar to  $f_3\wedge f_1$.

For $f_3\wedge f_3'$,  we have  $(f,g)_w=0$.

\ \

 Now, we check the compositions of left and right multiplication. We prove only the cases of  $
f_1P(u)$  and $P(u)f_3$, where $f_1\in F_1, \ f_3\in F_3$, $u\in
\Phi(X)$. Others can be similarly proved.

We may assume that $u\in \Phi_{1}(X)$.

For $f_1P(u)$, let $f=x_iP(x_j)-\{x_i\prec x_j\}, i,j\in I$ and
$w=\overline{fP(u)}$. There are two cases to consider.
\begin{enumerate}
\item[(I)]\ $u=x_k\in X$. Then by using the equation (\ref{e1}),
\begin{eqnarray*}
fP(u)&=&x_iP(x_j)P(x_k)-\{x_i\prec x_j\}P(x_k)\\
&=&x_iP(P(x_j)x_k)+x_iP(x_jP(x_k))-\{x_i\prec x_j\}P(x_k)\\
&\equiv&x_iP(\{x_j\succ x_k\})+x_iP(\{x_j\prec x_k\})-\{\{x_i\prec x_j\}\prec x_k\}\\
&\equiv&\{x_i\prec\{x_j\succ x_k\}\}+\{x_i\prec\{x_j\prec x_k\}\}-\{\{x_i\prec x_j\}\prec x_k\}\\
&\equiv&0\ \  mod(S, w).
\end{eqnarray*}
\item[(II)]\ $u\in \Phi_{1}(X)\setminus X$. Then
\begin{eqnarray*}
fP(u)&=&x_iP(x_j)P(u)-\{x_i\prec x_j\}P(u)\\
&=&x_iP(P(x_j)u)+x_iP(x_jP(u))-\{x_i\prec x_j\}P(u).
\end{eqnarray*}
By Lemma \ref{l2}, we have $P(x_j)u\equiv \Sigma\alpha_lu_l
 \ \ mod(F_{1}\cup F_{2}, \overline{P(x_j)u)}$, where $u_l\in \Phi_{1}(X)\setminus
 X$. Then
\begin{eqnarray*}
fP(u)&\equiv&x_iP(\Sigma \alpha_lu_l)+x_iP(x_jP(u))-\{x_i\prec x_j\}P(u)\\
&\equiv&\Sigma\alpha_lx_iP(u_l)+x_iP(x_jP(u))-\{x_i\prec x_j\}P(u)\\
&\equiv&0\ \  mod(S, w).
\end{eqnarray*}
\end{enumerate}

For $P(u)f_3$, let $f=P(v_1)u_{1}\cdots P(v_{n})u_{n}$, $u_{n}\in
X^*,\ u_{t}\in X^*\backslash \{1\},\ v_{l}\in\Phi_{1}(X)\setminus
X,\ n\geq1,\ 1\leq t<n,\ 1\leq l\leq n$, and let
$w=\overline{P(u)f}$. Then

$$
P(u)f=P(u)P(v_1)u_{1}\cdots P(v_{n})u_{n}.
$$

By Lemma \ref{l2}, we have $P(u)P(v_1)\equiv\Sigma\alpha_{i}P(w_{i})
\  mod(F_{1}\cup F_{2}, \overline{P(u)P(v_1)})$, where each
$w_{i}\in \Phi_{1}(X)\setminus X$. Then

\begin{eqnarray*}
P(u)f&\equiv&\Sigma\alpha_{i}P(w_{i})u_{1}\cdots
P(v_{n})u_{n}\\
&\equiv&0\ \  mod(S, w).
\end{eqnarray*}

So, all compositions in $S$ are trivial. The proof is complete. $
\hfill \blacksquare $

\ \

We reach to prove  Theorem \ref{t1}.

\ \

{\bf The proof of Theorem \ref{t1}:} Let $R=F_1\cup F_2$. Then for
any $u\not\in Irr(R^{comp})$, we have $u= v|_{\bar{r}}$, where $r\in
R^{comp}, v|_{r}\ \mbox{is normal}\ R^{comp}\mbox{-word}$.
 Then $f=v|_{r}\in Id(R^{comp})=Id(R)\subseteq Id(S)$. Since $S$ is a
Gr\"{o}bner-Shirshov basis in $RB(X)$, by Theorem \ref{CDL}, we have
$\bar{f}=w|_{\overline{s}}$ for some $w \in \Phi^{\star}(X)$,\ $s\in
S$. That is, $u= v|_{\bar{r}}=\bar{f}\not\in Irr(S)$. So, we have
that $Irr(R^{comp})\supset Irr(S)\supset X$. Since $Irr(R^{comp})$
is a $F$-basis of $U(D)$, $D$ can be embedded into $U(D)$.
 $\hfill \blacksquare $

\ \

\noindent{\bf Acknowledgement}: The authors would like to thank
Professor L.A. Bokut for his guidance, useful discussions and
enthusiastic encouragement given to this paper. The authors also
thank Professor L. Guo for his valuable suggestions to this paper.


\begin{thebibliography}{99}

\bibitem{AGKO}G.E. Andrew, L. Guo, W. Keigher, K. Ono, Baxter algebras and Hopf
algebras, {\it Trans. Amer. Math. Soc}., 355 (2003), 4639-4656.

\bibitem{A} M. Aymon and P.-P. Grivel, Un theoreme de Poincare-Birkhoff-Witt pour
les algebres de Leibniz, {\it Comm. Algebra},  31(2003), 527-544.

\bibitem{Baxter} G. Baxter, An analytic problem whose solution follows from a simple
algebraic identity, {\it Pacific J. Math.}, 10(1960), 731-742.

\bibitem{Birkhoff}G.D. Birkhoff, Representability of Lie algebras and Lie groups by
matrices, {\it Ann. of Math.}, 38(2)(1937), 526-532.

\bibitem{BCD}L.A. Bokut, Yuqun Chen and Xueming Deng, Gr\"{o}bner-Shirshov bases for Rota-Baxter
algebras, {\it Siberian Math. J.}, to appear.
arxiv.org/abs/0908.2281

\bibitem{Cartier} P. Cartier, On the structure of free Baxter algebras, {\it Adv. Math.}, 9(1972), 253-265.

\bibitem{Cohn}P.M. Cohn, A remark on the Birkhoff-Witt theorem, {\it J. London Math.
Soc.}, 38(1963), 197-203.

\bibitem{E02}K. Ebrahimi-Fard, Loday-type algebras and the
Rota-Baxter relation, {\it Lett. Math. Phys.}, 61(2)(2002), 139-147.

\bibitem{EG08a}K. Ebrahimi-Fard and L. Guo, Rota-Baxter algebras and dendriform algebras,
{\it Journal of Pure and Applied Algebra}, {212}(2)(2008), 320-339.

\bibitem{EG08b} K. Ebrahimi-Fard and L. Guo, Free Rota-Baxter algebras and rooted
trees,  {\it J. Algebra and Its Applications}, 7(2008), 167-194.



\bibitem{Higgins}P.J. Higgins, Baer Invariants and the Birkhoff-Witt theorem, {\it J. of
Algebra}, 11(1969), 469-482.


\bibitem{K} P.S. Kolesnikov, Conformal representations of Leibniz algebras,
{\it Siberian Math. J.}, 49(2008), 429-435.



\bibitem{Lo}J.-L. Loday, Dialgebras, in dialgebras and related
operads, {\it Lecture Notes in Math.}, {1763}(2001), 7-66.



\bibitem{Lo04}J.-L. Loday, Scindement d'associativite et algebres de Hopf. to appear in the
Proceedings of the Conference in honor of Jean Leray, Nantes(2002),
{\it Seminaire et Congres(SMF)}, 9(2004), 155-172.



\bibitem{Poincare}H. Poincare, Sur les groupes continus,  {\it Trans. Cambr. Philos.
Soc.}, 18(1900), 220-225.

\bibitem{Rota}G. Rota, Baxter operators, an introduction, In: ``Gian-Carlo Rota on
Combinatorics, Introduc- tory papers and commentaries", Joseph P.S.
Kung, Editor, Birkh¡§auser, Boston, 1995.

\bibitem{Shest} I.P. Shestakov, Every Akivis algebra is linear, {\it
Geometriae Dedicata}, 77(1999), 215-223.

\bibitem{Witt}E. Witt, Treue Darstellung Liescher Ringe, {\it J. Reine Angew.
Math.}, 177(1937), 152-160.

\end{thebibliography}
\end{document}